\documentclass[12pt]{article}
\usepackage[utf8]{inputenc}
\usepackage{amssymb,amsmath,amsfonts,eucal,mathrsfs,amsthm} 
\usepackage[colorinlistoftodos]{todonotes}
\usepackage[allcolors=blue]{hyperref}
\newtheorem{theorem}{Theorem}
\newtheorem{proposition}[theorem]{Proposition}

\newtheorem{corollary}[theorem]{Corollary}

\newtheorem{remark}[theorem]{Remark}

\theoremstyle{definition}
\newcommand{\R}{\mathbb{R}}

\newcommand{\spa}{\mbox{span}}

\newcommand{\Ima}{\mbox{Im }}

\newcommand{\nap}{\nabla^{\perp}}
\newcommand{\nab}{\tilde\nabla}

\newcommand{\End}{\mbox{End}}
\newcommand{\Hom}{\mbox{Hom}}

\newcommand{\Y}{\mathcal{Y}\,}
\newcommand{\D}{\mathcal{D}}
\def\<{{\langle}}
\def\>{{\rangle}}

\def\F{{\cal F}}

\def\T{{\cal T}}

\def\Y{{\cal Y}}
\def\n{\nabla}
\def\d{\partial}
\def\a{\alpha}

\def\be{\begin{equation} }
\def\ee{\end{equation} }
\def\Ima{Im}
\def\proof{\noindent{\it Proof:  }}
\def\qed{\ifhmode\unskip\nobreak\fi\ifmmode\ifinner
\else\hskip5 pt \fi\fi\hbox{\hskip5 pt \vrule width4 pt
height6 pt  depth1.5 pt \hskip 1pt }}

\makeatletter
\newcommand{\subjclass}[2][]{%
  \let\@oldtitle\@title%
  \gdef\@title{\@oldtitle\footnotetext{#1 \emph{Mathematics Subject Classification:} #2}}%
}
\newcommand{\keywords}[1]{%
  \let\@@oldtitle\@title%
  \gdef\@title{\@@oldtitle\footnotetext{\emph{Key words and phrases.} #1.}}%
}
\makeatother

\begin{document}

\title{Infinitesimal variations of submanifolds}
\author{M. Dajczer and M. I. Jimenez}
\date{}
\maketitle
\begin{abstract} This paper deals with the subject of infinitesimal 
variations of Euclidean submanifolds with arbitrary dimension and 
codimension. 
The main goal is to establish a Fundamental theorem for these 
geometric objects. Similar to the theory of isometric immersions 
in Euclidean space,  we prove that a system of three equations 
for a certain pair of tensors  are the integrability conditions 
for the differential equation that determines the infinitesimal 
variations.
In addition, we give some rigidity results when the submanifold 
is intrinsically a Riemannian product of manifolds.
\end{abstract}

The classical theory of smooth variations of surfaces in $\R^3$ was 
the object of intense study by geometers in the $19^{th}$ century,
initially with no distinction between isometric and infinitesimally 
isometric variations. But by the end of that century a clear 
distinction was done by Darboux.  In the $20^{th}$ century a vast 
number of papers on the subject was produced, as can be seen in the 
survey paper \cite{IMS}. In particular, it is quoted that Efimov 
observed that  ``the theory of infinitesimal variations is the 
differential of the theory of isometric variations". For a modern 
account of several aspects of the subject we refer to Spivak \cite{Sp}.   

The case of infinitesimal variations of Euclidean hypersurfaces 
has been initially considered around the beginning of the $20^{th}$ 
century, and the work of Sbrana \cite{Sb} in 1908 stands out. 
A complete local parametric classification of the hypersurfaces 
that admit infinitesimal variations is due to Dajczer and Vlachos 
\cite{DV} 
who, in particular, showed that this class is much larger than 
the one of hypersurfaces that allow isometric variations. The 
classification in the case of complete hypersurfaces 
is due to Jimenez \cite{Ji}. 

Dajczer and Rodríguez \cite{DR} showed that Euclidean submanifolds 
in low codimension are generically infinitesimally rigid, that is, 
only trivial infinitesimal variations are possible. An infinitesimal 
variation is called trivial if it relates to a variation that is 
the restriction to the submanifold of a smooth one-parameter family 
of isometries of 
Euclidean space.  For hypersurfaces, their result is already 
contained in the book of Cesàro \cite{Ce} published in 1886.
For recent results on the subject of rigidity of infinitesimal
variations we refer to Dajczer and Jimenez \cite{DJ}.

Our goal in this paper is twofold.  In the first part, we obtain
a clean version of the Fundamental theorem of infinitesimal 
variations that extends to any codimension the result for 
hypersurfaces given in \cite{DV}.  A quite cumbersome coordinate 
version  for general codimension has been stated in \cite{Ma}. 
Similarly, as in the theory of isometric immersions, in terms of 
a pair of tensors associated to the variation, we obtain a system 
of three equations that are shown to be the integrability 
conditions for the equations that determine the infinitesimal 
variations.

Although the Fundamental theorem is stated and proved for the
ambient Euclidean space with the standard Riemannian metric, the 
same  statement and proof work if the ambient Euclidean space has 
any possible signature. 

The second part of the paper is devoted to some rigidity results 
concerning the situation when the starting manifold is a Riemannian 
product of manifolds. We provide conditions that imply that any 
infinitesimal variation is given by infinitesimal variations of 
each factor.

\section{Preliminaries}

Let $f\colon M^n\to\R^m$ of an isometric immersion of an 
$n$-dimensional Riemannian manifold $M^n$ into Euclidean space with 
codimension $m-n$. An isometric  variation of $f$ is a  
smooth map $\F\colon I\times M^n\to\R^m$ for some open interval 
$0\in I\subset\R$ such that $f(x)=\F(0,x)$ and 
$f_t=\F(t,\cdot)\colon M^n\to\R^m$ is an isometric immersion for 
any $t\in I$. The variation is said 
to be trivial if it is produced by a smooth family of isometries 
of $\R^m$. That is, if there exist a smooth family 
$C\colon I\to O(m)$ of orthogonal transformations of $\R^m$ with 
$C(0)=I$ and a smooth map $v\colon I\to\R^m$ with $v(0)=0$, and 
$$
\F(t,x)=C(t)f(x)+v(t)
$$
for all $(t,x)\in I\times M^n$. 

The variational vector field of the isometric variation $\F$ is 
the section $\T\in\Gamma(f^*T\R^m)$ defined as 
$\T=\F_*\d/\d t|_{t=0}$. It is easily seen that it  satisfies 
the condition
\be\label{cis}
\<\nab_X\T,f_*Y\>+\<f_*X,\nab_Y\T\>=0
\ee
for any tangent vector fields $X,Y\in\mathfrak{X}(M)$. Here and in 
the sequel, we use the same notation for the inner products in $M^n$ 
and $\R^m$ and denote by $\n$ and $\nab$ the respective associated 
Levi-Civita connections.
\vspace{1ex}

In this paper, we deal with the concept of infinitesimal variation, 
that is, the infinitesimal analogue to an isometric  variation.
These variations preserve lengths just ``up to the first order".
\vspace{1ex}

A smooth variation $F\colon I\times M^n\to\R^m$ of a given isometric 
immersion $f\colon M^n\to\R^m$ is called an \emph{infinitesimal variation} 
if it satisfies the condition
\be\label{varcond}
\frac{\d}{\d t}|_{t=0}\<f_{t*}X,f_{t*}Y\>=0
\ee
for any $X,Y\in\mathfrak{X}(M)$. 
\vspace{1ex}

In order to study the infinitesimal analogue of isometric variations 
it is known from classical geometry that the convenient approach 
is to look at the variational vector field. In our 
case, that this is the way to proceed is justified in the sequel.
\vspace{1ex}

We have from \eqref{varcond} that the variational vector field $\T$ 
of $F$ has to satisfy condition \eqref{cis}. This leads to the 
following definition.
\vspace{1ex}

A section $\T$ of $f^*T\R^m$ is called an \emph{infinitesimal 
bending} of an isometric immersion $f\colon M^n\to\R^m$ if 
the condition
\be\label{cib}
\<\nab_X\T,f_*Y\>+\<f_*X,\nab_Y\T\>=0
\ee
holds for any tangent vector fields $X,Y\in\mathfrak{X}(M)$.
\vspace{1ex}

Associated to an infinitesimal bending $\T$ of $f$ we have that 
the infinitesimal variation 
$\F\colon\R\times M^n\to\R^m$ given by
\be\label{unique}
\F(t,x)=f(x)+t\T(x)
\ee
has variational vector field $\T$. The statement that the 
length is preserved up to the first order is because
$$
\|f_{t*}X\|^2=\|f_*X\|^2+t^2\|\tau_*X\|^2
$$
for any $X\in\mathfrak{X}(M)$. By no means
\eqref{unique} is unique with this property, although it may 
be seen as the simplest one.  In fact, new  
infinitesimal variations with variational vector field $\T$ 
are obtained by adding to \eqref{unique} terms of the type 
$t^k\delta$, $k>1$, where $\delta\in\Gamma(f^*T\R^{n+1})$ and, 
maybe, for restricted values of the parameter $t$.
\vspace{1ex}

A \emph{trivial infinitesimal bending} is the restriction 
to the submanifold of a Killing vector field of the ambient space.
More precisely,  there is a skew-symmetric linear 
endomorphism ${\cal D}$ of $\R^m$ and $v\in\R^m$ such that
${\cal T}={\cal D}f+v$. Then, we have that
$$
\F(t,x)=e^{t{\cal D}}f(x)+tv
$$
is a trivial isometric variation of $f$. Throughout the paper
we identify pairs of infinitesimal bendings that differ by
a trivial one.

\section{The Fundamental theorem}

In this section, we first define a pair of tensors associated 
to any infinitesimal bending and then show that they satisfy 
a set of three equations that form the \emph{Fundamental system 
of equations} of the bending. Then we state and prove the 
\emph{Fundamental theorem of infinitesimal variations} by showing 
that the equations of the system are the integrability conditions 
for the system of equations whose solutions give infinitesimal 
bendings. 
\vspace{1ex}

Let $\T\in\Gamma(f^*T\R^m)$ be an  infinitesimal bending of a given 
isometric immersion $f\colon M^n\to\R^m$. We first argue that $\T$ 
and the second fundamental form $\a\colon TM\times TM\to N_fM$ 
of $f$ determine an \emph{associate pair} of tensors 
$(\beta,{\cal E})$, where $\beta\colon TM\times TM\to N_fM$ 
is symmetric and ${\cal E}\colon TM\times N_fM\to N_fM$ 
satisfies the compatibility condition
\be\label{anti}
\<{\cal E}(X,\eta),\xi\>+\<{\cal E}(X,\xi),\eta\>=0
\ee
for any $X\in\mathfrak{X}(M)$ and $\eta,\xi\in\Gamma(N_fM)$.
\vspace{1ex}

Let $L\in\Gamma(\Hom(TM,f^*T\R^m))$ be the tensor defined by
$$
LX=\nab_X \T=\T_*X
$$
for any $X\in\mathfrak{X}(M)$.
Notice that \eqref{cib} in terms  of $L$ has the form
\be\label{inf}
\<LX,f_*Y\>+\<f_*X,LY\>=0
\ee
for any $X,Y\in\mathfrak{X}(M)$. Let 
$B\colon TM\times TM\to f^*T\R^m$ be the tensor given by
$$
B(X,Y)=(\nab_XL)Y=\nab_XLY-L\n_XY
$$
for any $X,Y\in\mathfrak{X}(M)$. Since 
$$
B(X,Y)=\nab_X\nab_Y\T-\nab_{\nabla_XY}\T
$$
then $B$ is symmetric due to the  flatness of the ambient space.
Then also the tensor $\beta\colon TM\times TM\to N_fM$  defined by
$$
\beta(X,Y)=(B(X,Y))_{N_fM}
$$
is symmetric. 

For later use, we define the symmetric tensor 
$B_\xi\in\Gamma(\End(TM))$  by
$$
\<B_\xi X,Y\>=\<\beta(X,Y),\xi\>
$$
for any $X,Y\in\mathfrak{X}(M)$ and $\xi\in\Gamma(N_fM)$.

Let $\Y\in\Gamma(\Hom(N_fM,f_*TM))$ be given by
\be\label{Y}
\<\Y\eta,f_*X\>+\<\eta,LX\>=0.
\ee
Then, the tensor ${\cal E}\colon TM\times N_fM\to N_fM$ 
is defined by
$$
{\cal E}(X,\eta)=\a(X,\Y\eta)+(LA_\eta X)_{N_fM}
$$
where $A_\eta\in\Gamma(\End(TM))$ is given by
$$
\<A_\eta X,Y\>=\<\a(X,Y),\eta\>.
$$
We have
\begin{align*}
\<{\cal E}(X,\eta),\xi\>
&=\<\a(X,\Y\eta)+LA_\eta X,\xi\>
=\<A_\xi X,\Y\eta\>-\<\Y\xi,A_\eta X\>\\
&=-\<LA_\xi X,\eta\>-\<\a(X,\Y\xi),\eta\>
=-\<{\cal E}(X,\xi),\eta\>,
\end{align*}
and hence \eqref{anti} is satisfied.

\begin{proposition} We have that
\be\label{tam}
(B(X,Y))_{TM}=\Y\a(X,Y)
\ee
for any $X,Y\in\mathfrak{X}(M)$.
\end{proposition}

\proof We need to show that
$$
C(X,Y,Z)=\<(B-\Y\a)(X,Y),f_*Z\>
$$
vanishes for any $X,Y,Z\in\mathfrak{X}(M)$. The 
derivative of \eqref{inf} gives
\begin{align*}
0&=\<\nab_ZLX,f_*Y\>+\<LX,\nab_Zf_*Y\>
+\<\nab_ZLY,f_*X\>+\<LY,\nab_Zf_*X\>\\ 
&=\<B(Z,X),f_*Y\>+\<L\n_ZX,f_*Y\>+\<LX,f_*\n_ZY+\a(Z,Y)\>\\
&\;\;\;+\<B(Z,Y),f_*X\>+\<L\n_ZY,f_*X\>+\<LY,f_*\n_ZX+\a(Z,X)\>\\
&=\<B(Z,X),f_*Y\>+\<LX,\a(Z,Y)\>+\<B(Z,Y),f_*X\>+\<LY,\a(Z,X)\>\\
&=\<(B-\Y\a)(Z,X),f_*Y\>+\<(B-\Y\a)(Z,Y),f_*X\>.
\end{align*}
From the symmetry of $B$ and the above, we obtain 
$$
C(X,Y,Z)=C(Y,X,Z)\;\;\;\text{and}\;\;\;C(Z,X,Y)=-C(Z,Y,X)
$$
for any $X,Y,Z\in\mathfrak{X}(M)$. Then
\begin{align*}
C(X,Y,Z)&=-C(X,Z,Y)=-C(Z,X,Y)=C(Z,Y,X)\\
&=C(Y,Z,X)=-C(Y,X,Z)=-C(X,Y,Z)=0,
\end{align*}
as we wished.\vspace{2ex}\qed

The last manipulation of the above proof is known as
the Braid Lemma; for instance see \cite{Be} p.\ $224$.

\begin{proposition}\label{system}
The pair of tensors $(\beta,{\cal E})$ associated to an 
infinitesimal bending $\T$ satisfies the following system 
of three equations:
\be\label{Gauss}
A_{\beta(Y,Z)}X+B_{\alpha(Y,Z)}X=A_{\beta(X,Z)}Y
+B_{\alpha(X,Z)}Y,
\ee
\be\label{Codazzi}
(\nabla^{\perp}_X\beta)(Y,Z)-(\nabla_Y^{\perp}\beta)(X,Z)
={\cal E}(Y,\a(X,Z))-{\cal E}(X,\a(Y,Z))
\ee
and
\begin{align}\label{Ricci}
(\nap_X{\cal E})(Y,\eta)&-(\nap_Y{\cal E})(X,\eta)\nonumber\\
&=\beta(X,A_\eta Y)-\beta(A_\eta X,Y)+\a(X,B_\eta Y)-\a(B_\eta X,Y)
\end{align}
for all $X,Y,Z\in\mathfrak{X}(M)$ and $\eta\in\Gamma(N_fM)$.
Moreover \eqref{Codazzi} is equivalent to
\be\label{Codazzi2}
(\nabla_XB_\eta)Y-(\nabla_YB_\eta)X-B_{\nap_X\eta}Y
+B_{\nap_Y\eta}X=A_{{\cal E}(X,\eta)}Y-A_{{\cal E}(Y,\eta)}X
\ee
for all $X,Y,Z\in\mathfrak{X}(M)$ and $\eta\in\Gamma(N_fM)$.
\end{proposition}

\proof  We first show that
\be\label{derY}
(\nab_X\Y)\eta=-f_*B_\eta X-LA_\eta X+{\cal E}(X,\eta)
\ee
for any $X\in\mathfrak{X}(M)$ and $\eta\in\Gamma(N_fM)$,
where $(\nab_X \Y)\eta=\nab_X\Y\eta-\Y\nap_X\eta$. 
We have from \eqref{Y} and the derivative of \eqref{inf} that
\begin{align*}
0&=\<\nab_X\Y\eta,f_*Y\>+\<\Y\eta,f_*\n_XY\>
+\<\nab_XLY,\eta\>+\<LY,\nab_X\eta\>\\
&=\<(\nab_X\Y)\eta,f_*Y\>+\<B_\eta X,Y\>+\<LA_\eta X,f_*Y\>.
\end{align*} 
Since $\<\Y\eta,\xi\>=0$, we obtain
\begin{align*}
 0&=\<\nab_X\Y\eta,\xi\>+\<\Y\eta,\nab_X\xi\>
 =\<(\nab_X\Y)\eta,\xi\>-\<\a(X,\Y\eta),\xi\>\\
 &=\<(\nab_X\Y)\eta,\xi\>+\<LA_\eta X-{\cal E}(X,\eta),\xi\>
\end{align*}
for any $X\in\mathfrak{X}(M)$ and $\eta,\xi\in\Gamma(N_fM)$, 
and \eqref{derY} follows.

Using
\be\label{form}
(\nab_XB)(Y,Z)
=\nab_X(\nab_YL)Z-(\nab_{\nabla_XY}L)Z-(\nab_Y L)\nabla_XZ
\ee
it is easy to see that 
\be\label{segderL}
(\nab_XB)(Y,Z)-(\nab_YB)(X,Z)=-LR(X,Y)Z
\ee
for all $X,Y,Z\in\mathfrak{X}(M)$. 
It follows using \eqref{tam} that
$$
\<(\nab_X B)(Y,Z),f_*W\>=\<(\nab_X \Y)\a(Y,Z)+\Y(\nap_X \a)(YZ)
-f_*A_{\beta(Y,Z)}X,f_*W\>
$$
for any $X,Y,Z,W\in\mathfrak{X}(M)$. Then
\eqref{segderL} and the Codazzi equation give
\begin{align*}
\<(\nab_X \Y)\a(Y,Z)&-(\nab_Y \Y)\a(X,Z),f_*W\>\\
&=\<LR(Y,X)Z+A_{\beta(Y,Z)}X-A_{\beta(X,Z)}Y,W\>
\end{align*}
Now using  the Gauss equation, we obtain
\begin{align*}
\<(\nab_X \Y)\a(Y,Z)&-(\nab_Y \Y)\a(X,Z),f_*W\>\\
&=\<LA_{\a(X,Z)}Y-LA_{\a(Y,Z)}X
+A_{\beta(Y,Z)}X-A_{\beta(X,Z)}Y,f_*W\>.
\end{align*}
On the other hand, it follows from \eqref{derY} that
\begin{align*}
\<(\nab_X \Y)\a(Y,Z)&-(\nab_Y \Y)\a(X,Z),f_*W\>\\
&=\<B_{\a(X,Z)}Y+LA_{\a(X,Z)}Y-B_{\a(Y,Z)}X-LA_{\a(Y,Z)}X,f_*W\>.
\end{align*}
From the last two equations, we obtain
\begin{align*}
\<B_{\a(X,Z)}Y&-B_{\a(Y,Z)}X,f_*W\>=
\<A_{\beta(Y,Z)}X-A_{\beta(X,Z)}Y,W\>,    
\end{align*}
and this is \eqref{Gauss}.

Using \eqref{form} we obtain
$$
((\nab_XB)(Y,Z))_{N_fM}=
\alpha(X,\Y\a(Y,Z))+(\nabla_X^{\perp}\beta)(Y,Z).
$$
Then, we have from \eqref{segderL} and the Gauss equation that
\begin{align*}
(\nabla_X^{\perp}&\beta)(Y,Z)-(\nabla_Y^{\perp}\beta)(X,Z)\\
&=(LR(Y,X)Z)_{N_fM}-\alpha(X,\Y\a(Y,Z)+\alpha(Y,\Y\a(X,Z)\\
&=(LA_{\a(X,Z)}Y-LA_{\a(Y,Z)}X)_{N_fM}
-\alpha(X,\Y\a(Y,Z)+\alpha(Y,\Y\a(X,Z),
\end{align*}
and this is \eqref{Codazzi}. Since ${\cal E}$ satisfies the 
compatibility condition \eqref{anti}, then
$$
\<{\cal E}(X,\a(Y,Z)),\eta\>=-\<A_{{\cal E}(X,\eta)}Y,Z\>.
$$
and this gives \eqref{Codazzi2}.

We have
\begin{align*}
(\nap_X{\cal E})(Y,\eta)&=\nap_X{\cal E}(Y,\eta)
-{\cal E}(\nabla_XY,\eta)-{\cal E}(Y,\nap_X\eta)\\
&=(\nap_X\a)(Y,\Y\eta)+(L(\nabla_X A)(Y,\eta))_{N_fM}
+\a(Y,\nabla_X\Y\eta)\\
&\;\;\;-\a(Y,\Y\nap_X\eta)-(L\nabla_XA_\eta Y)_{N_fM}
+\nap_X(LA_\eta Y)_{N_fM}.
\end{align*}
Then \eqref{derY} gives
\begin{align*}
(\nap_X{\cal E})(Y,\eta)&=(\nap_X\a)(Y,\Y\eta)
+(L(\nabla_X A)(Y,\eta))_{N_fM}-\a(Y,B_\eta X)\\
&\;\;\;-\a(Y,(LA_\eta X)_{TM})-(L\nabla_XA_\eta Y)_{N_fM}
+\nap_X(LA_\eta Y)_{N_fM}.
\end{align*}

Using the Codazzi equation, we obtain
\begin{align*} 
(\nap_X{\cal E})(Y,\eta)&-(\nap_Y{\cal E})(X,\eta)=
\a(X,B_\eta Y)-\a(Y,B_\eta X)+\a(X,(LA_\eta Y)_{TM})\\
&-\a(Y,(LA_\eta X)_{TM})-(L\nabla_XA_\eta Y)_{N_fM}
+\nap_X(LA_\eta Y)_{N_fM}\\
&+(L\nabla_YA_\eta X)_{N_fM}-\nap_Y(LA_\eta X)_{N_fM}.
\end{align*}
Since
$$
\beta(X, A_\eta Y)=\a(X,(LA_\eta Y)_{TM})-(L\nabla_XA_\eta Y)_{N_fM}
+\nap_X(LA_\eta Y)_{N_fM},
$$
then \eqref{Ricci} follows.\vspace{2ex}\qed

We say that an isometric immersion $f\colon M^n \to\R^m$ has
\emph{full}  first normal space $N_1^f(x)$ at $x\in M^n$ if 
$$
N_1^f(x)=\spa\{\a(X,Y): X,Y\in T_xM\}
$$
satisfies $N_1^f(x)=N_fM(x)$.
\vspace{1ex}

The following result shows that for a submanifold with full 
first normal spaces the tensor $\beta$ determines $\mathcal{E}$.

\begin{proposition}\label{uniqueE}
Let $f\colon M^n\to\R^m$ be an isometric immersion with full 
first normal spaces. If $(\beta, \mathcal{E})$ is the associated 
pair of tensors to an infinitesimal bending $\T$ of $f$ then 
$\mathcal{E}$ is the unique tensor that satisfies \eqref{anti} 
and \eqref{Codazzi}.
\end{proposition}

\proof If $\mathcal{E}_0\colon TM\times N_fM\to N_fM$ is a tensor 
that satisfies \eqref{anti} and \eqref{Codazzi}, it follows from 
\eqref{Codazzi} that
$$
(\mathcal{E}-\mathcal{E}_0)(X,\a(Y,Z))
=(\mathcal{E}-\mathcal{E}_0)(Y,\a(X,Z))
$$
for any $X,Y,Z\in\mathfrak{X}(M)$. Since both $\mathcal{E}$
and $\mathcal{E}_0$ satisfy \eqref{anti}, we have 
$$
\<(\mathcal{E}-\mathcal{E}_0)(X_1,\a(X_2,X_3)),\a(X_4,X_5)\>
\!=\!-\<(\mathcal{E}-\mathcal{E}_0)(X_1,\a(X_4,X_5)),\a(X_2,X_3)\>
$$
where $X_i\in\mathfrak{X}(M)$, $1\leq i\leq 5$. We denote
$$
\<(\mathcal{E}-\mathcal{E}_0)(X_1,\a(X_2,X_3)),\a(X_4,X_5)\>
=(X_1,X_2,X_3,X_4,X_5).
$$
It follows from the relations above and the symmetry of $\a$ 
that
\begin{align*}
&(X_1,X_2,X_3,X_4,X_5)=-(X_1,X_4,X_5,X_2,X_3)=-(X_5,X_4,X_1,X_3,X_3)\\
&=(X_5,X_2,X_3,X_4,X_1)=(X_3,X_2,X_5,X_4,X_1)=-(X_3,X_4,X_1,X_2,X_5)\\
&=-(X_4,X_3,X_1,X_2,X_5)=(X_4,X_2,X_5,X_3,X_1)=(X_2,X_4,X_5,X_3,X_1)\\
&=-(X_2,X_3,X_1,X_4,X_5)=-(X_2,X_1,X_3,X_4,X_5)=-(X_1,X_2,X_3,X_4,X_5)\\
&=0,
\end{align*}
and thus $\mathcal{E}-\mathcal{E}_0=0$.\qed

\begin{remark}{\em 
An alternative way to obtain the equations in Proposition
\ref{system} is to follow the ``classical" procedure,
which goes as follows.  Since the metrics $g_t$ induced by
the infinitesimal variation $f_t=f+t\T$ satisfy
$\d/\d t|_{t=0}g_t=0$, hence the Levi-Civita
connections and curvature tensors of $g_t$ satisfy
$$
\d/\d t|_{t=0}\nabla^{t}_X Y=0
$$
and
$$
\d/\d t|_{t=0}g_t(R^t(X,Y)Z,W)=0
$$
for any $X,Y,Z,W\in\mathfrak{X}(M)$. Then use this to compute 
the derivatives with respect to $t$ at $t=0$ of the Gauss, 
Codazzi and Ricci equations for $f_t$. In fact, as can be seen 
in \cite{DJ} this works quite nicely to obtain \eqref{Gauss}. 
On the contrary, the computation for the other two equations becomes 
really cumbersome outside the hypersurface case. For hypersurfaces
this was done in \cite{DT} and \cite{DV}. A result in coordinates  
for general codimension has been stated in \cite{Ma}.
}\end{remark}

Let $f\colon M^n\to\R^m$ be an isometric immersion and $\mathcal{T}$ 
be a trivial infinitesimal bending of $f$, that is,
$$
\T=\D f+w
$$
where $\D\in\End(\R^m)$ is skew-symmetric and $w\in\R^m$. Then
$$
L=\D|_{f_*TM}\;\;\text{and}\;\;B(X,Y)=\D\a(X,Y).
$$
Let $\D^N\in\Gamma(\End(N_fM))$ skew-symmetric be given by
$$
\D^N\eta=(\D\eta)_{N_fM}
$$
for any $\eta\in\Gamma(N_fM)$. Then we have
$$
\beta(X,Y)=\D^N\a(X,Y)\;\;\text{and}\;\;{\cal E}(X,\eta)
=-(\nap_X\D^N)\eta,
$$  
where the second equation follows computing $(\nab_XD)\eta=0$.

\begin{proposition}\label{trivial}
An infinitesimal bending $\mathcal{T}$ is trivial if and only 
if there is $C\in\Gamma(\End(N_fM))$ skew-symmetric such that 
\be\label{contriv}
\beta(X,Y)=C\a(X,Y)\;\;\text{and}\;\;{\cal E}(X,\eta)
=-(\nap_XC)\eta.
\ee 
\end{proposition}

\proof Define  $\D\in\Gamma(\End(f^*\T\R^m))$ by
$$
\D(x)X=L(x)X\;\;\text{and}\;\;
\D(x)\eta=\Y(x)\eta+C(x)\eta
$$
for any $X\in T_xM$ and $\eta\in N_{f(x)}M$.
Using the assumption on $\beta$ we obtain
\begin{align*}
\nab_X\D Y&=(\nab_X L)Y+L\nabla_XY
=\Y\a(X,Y)+C\a(X,Y)+L\nabla_XY\\
&=\D\nab_XY
\end{align*}
for any $X,Y\in\mathfrak{X}(M)$. 
The assumptions on ${\cal E}$ and \eqref{derY} give
\begin{align*}
\nab_X\D\eta&=\nab_X\Y\eta+\nab_XC\eta\\
&=(\nab_X\Y)\eta+\Y\nap_X\eta+(\nap_XC)\eta
+C\nap_X\eta-f_*A_{C\eta}X\\
&=-f_*B_\eta X-LA_\eta X+\Y\nap_X\eta
+C\nap_X\eta-f_*A_{C\eta}X
\end{align*}
for any $X\in\mathfrak{X}(M)$ and $\eta\in\Gamma(N_fM)$.
But $B_\eta=-A_{C\eta}$ from $\beta=C\a$, then
$$
\nab_X\D\eta=-LA_\eta X+\Y\nap_X\eta+C\nap_X\eta
=\D\nab_X\eta.
$$
Therefore, we have shown that $\D(x)=\D$ is constant along 
$M^n$. Thus the map $\mathcal{T}-\D f$ is constant.
\vspace{2ex}\qed

Recall that we identify two infinitesimal bendings $\T_1$ and 
$\T_2$ of $f$ whenever $\T_0=\T_2-\T_1$ is a trivial infinitesimal 
bending. 
In this case, if $(\beta_1,{\cal E}_1)$ and $(\beta_2,{\cal E}_2)$ 
are the associated pairs to $\T_1$ and $\T_2$, respectively, 
then the associated pair to $\T_0$ is 
$(\beta_2-\beta_1,{\cal E}_2-{\cal E}_1)$ that
satisfies \eqref{contriv} for $\D^N\in\Gamma(\End(N_fM))$ 
skew-symmetric. Thus, in this situation we identify 
$(\beta_1,{\cal E}_1)$ and $(\beta_2,{\cal E}_2)$.
\vspace{1ex}

The following is the Fundamental theorem of infinitesimal variations.

\begin{theorem}\label{fund}Let $f\colon M^n\to\R^m$ be an 
isometric immersion of a simply connected Riemannian manifold. 
Let $\beta\colon TM\times TM\to N_fM$ be a symmetric tensor
and let the tensor ${\cal E}\colon TM\times N_fM\to N_fM$ 
satisfy the compatibility condition \eqref{anti}. 
If the pair $0\neq (\beta,{\cal E})$ satisfies 
\eqref{Gauss}, \eqref{Codazzi} and \eqref{Ricci}, 
then there is a unique infinitesimal bending $\mathcal{T}$ 
of $f$ having $(\beta,{\cal E})$ as associated pair. 
\end{theorem}

\proof Given a pair $(\beta,{\cal E})$ as in the statement, 
we first argue that there is $\D\in\Gamma(\End(f^*T\R^m))$ 
satisfying
\be\label{I}
(\nab_X \D)(Y+\eta)=-f_*B_\eta X+\beta(X,Y)+{\cal E}(X,\eta)
\ee
for any $X,Y\in\mathfrak{X}(M)$ and $\eta\in\Gamma(N_fM)$. 
To prove this, henceforth we check that the integrability 
condition of \eqref{I}, namely, that
$$
(\nab_X\nab_Y\D-\nab_Y\nab_X\D-\nab_{[X,Y]}\D)(Z+\eta)=0
$$
holds for any $X,Y,Z\in\mathfrak{X}(M)$ and $\eta\in\Gamma(N_fM)$. 
For simplicity, in the following we write $X$ instead of 
$f_*X$. We have 
\begin{align*}
&(\nab_X\nab_Y\D-\nab_Y\nab_X\D-\nab_{[X,Y]}\D)(Z+\eta)\\
&=\nab_X(\nab_Y\D)(Z+\eta)-(\nab_Y\D)\nab_X(Z+\eta)
-\nab_Y(\nab_X\D)(Z+\eta)\\
&\;\;\;+(\nab_X\D)\nab_Y(Z+\eta)-(\nab_{[X,Y]}\D)(Z+\eta)\\
&=\nab_X[-B_\eta Y+\beta(Y,Z)+\mathcal{E}(Y,\eta)]
+B_{\a(X,Z)+\nap_X\eta}Y-\beta(Y,\nabla_XZ-A_\eta X)\\
&\;\;\;-\mathcal{E}(Y,\a(X,Z)+\nap_X\eta)
+\nab_Y[B_\eta X-\beta(X,Z)-\mathcal E(X,\eta)]\\
&\;\;\;-B_{\a(Y,Z)+\nap_Y\eta}X
+\beta(X,\nabla_YZ-A_\eta Y)+\mathcal{E}(X,\a(Y,Z)+\nap_Y\eta)\\
&\;\;\;+B_\eta[X,Y]-\beta([X,Y],Z)-\mathcal{E}([X,Y],\eta).
\end{align*}
Hence
\begin{align*}
(\nab_X&\nab_Y\D-\nab_Y\nab_X\D-\nab_{[X,Y]}\D)(Z+\eta)\\
&=-A_{\beta(Y,Z)}X+B_{\a(X,Z)}Y+A_{\beta(X,Z)}Y-B_{\a(Y,Z)}X\\
&\;\;\;+(\nap_X\beta)(Y,Z)-(\nap_Y\beta)(X,Z)+\mathcal{E}(X,\a(Y,Z))
-\mathcal{E}(Y,\a(X,Z))\\
&\;\;\;-(\nabla_X B_\eta)Y+(\nabla_Y B_\eta)X+B_{\nap_X\eta}Y
-B_{\nap_Y\eta}X-A_{\mathcal{E}(Y,\eta)}X+A_{\mathcal{E}(X,\eta)}Y\\
&\;\;\;+(\nap_X\mathcal{E})(Y,\eta)-(\nap_Y\mathcal{E})(X,\eta)
-\a(X,B_\eta Y)+\a(Y,B_\eta X)\\
&\;\;\;+\beta(Y,A_\eta X)-\beta(X,A_\eta Y)\\
&=0,
\end{align*}
where for the last equality we made use \eqref{Gauss}, 
\eqref{Codazzi}, \eqref{Ricci}  and \eqref{Codazzi2}. 
\vspace{1ex}

Fix a solution $\D^*\in\Gamma(\End(f^*T\R^m))$ of \eqref{I}
and a point $x_0\in M^n$. Set $\D_0=\D^*(x_0)$ and let
\mbox{$\phi\colon f^*T\R^m\times f^*T\R^m\to\R$} be the tensor
defined by
$$
\phi(\rho,\sigma)
=\<(\D^*-\D_0)\rho,\sigma\>+\<(\D^*-\D_0)\sigma,\rho\>.
$$
Using \eqref{anti} and \eqref{I} we have
$(\nab_X\phi)(\rho,\sigma)=0$.
Hence $\phi=0$, and thus the maps $\D(x)=\D^*(x)-\D_0$
are skew-symmetric endomorphisms of $\R^m$.

Define  $L\in\Gamma(\Hom(TM,f^*T\R^m))$ by $L(x)=\D(x)|_{T_xM}$.
Using \eqref{I} we obtain
$$
(\nab_XL)Y=\nab_X\D Y-\D\nabla_XY=\beta(X,Y)+\D\a(X,Y).
$$
Then
$$
(\nab_XL)Y=(\nab_YL)X.
$$
Hence, there is $\mathcal{T}\in\Gamma(f^*T\R^m)$ such that
$$
\nab_X\mathcal{T}=LX
$$
for any $X\in\mathfrak{X}(M)$. Since $\D$ is
skew-symmetric then $L$ satisfies
$$
\<LX,Y\>+\<LY,X\>=0,
$$
proving that $\mathcal{T}$ is an infinitesimal
bending of $f$.  Moreover, its associate pair of tensors
$(\tilde{\beta}, \tilde{\mathcal{E}})$ is
$$
\tilde{\beta}(X,Y)=\beta(X,Y)+\D^N\a(X,Y)
\;\;
and\;\;
\tilde{\mathcal{E}}(X,\eta)=\mathcal{E}(X,\eta)
-(\nap_X\D^N)\eta.
$$
In fact, in this case $\Y\eta=(\D\eta)_{TM}$.
Using \eqref{I}, we have
\begin{align*}
\tilde{\mathcal{E}}(X,\eta)&=\a(X,(\D\eta)_{TM})
+(LA_{\eta}X)_{N_fM}\\
&=(\nab_X(\D\eta)_{TM})_{N_fM}+(LA_{\eta}X)_{N_fM}\\
&= (\nab_{X}\D\eta)_{N_fM}-\nap_X\D^N\eta
+(LA_{\eta}X)_{N_fM}\\
&=\mathcal{E}(X,\eta)+(\D\nab_X\eta)_{N_fM}
-\nap_X\D^N\eta+(LA_{\eta}X)_{N_fM}\\
&=\mathcal{E}(X,\eta)-(LA_\eta X)_{N_fM}
-(\nap_X\D^N)\eta+(LA_{\eta}X)_{N_fM}\\
&=\mathcal{E}(X,\eta)-(\nap_X\D^N)\eta.
\end{align*}
Another solution $\D^*_1$ of \eqref{I} gives rise 
to an infinitesimal bending $\T_1$ of $f$. It follows 
from Proposition \ref{trivial} that $\T-\T_1$ is a 
trivial infinitesimal bending.
\vspace{2ex}\qed

The following result is Theorem $13$ in \cite{DV}
or Theorem $14.11$ in \cite{DT}.

\begin{corollary}
Let $f\colon M^n\to\R^{n+1}$ be an isometric immersion of 
a simply connected Riemannian manifold. Let 
$0\neq {\cal B}\in\Gamma(\End(TM))$ be a symmetric 
Codazzi tensor that satisfies
$$
{\cal B}X\wedge AY-{\cal B}Y\wedge AX=0
$$
for all $X,Y\in\mathfrak{X}(M)$. Then there exists a unique 
infinitesimal bending $\T$ of $f$ having ${\cal B}$ 
as associated tensor.
\end{corollary}

\proof In this case the tensor $\mathcal{E}$ vanishes. Let 
$\beta\colon TM\times TM\to N_fM$ be the symmetric tensor 
given by $\beta(X,Y)=\<\mathcal{B}X,Y\>N$. Then \eqref{Ricci} 
trivially holds for $\beta$ and $\mathcal{E}=0$. Moreover, 
by the  assumptions on $\mathcal{B}$ we have that $(\beta,0)$ 
satisfies \eqref{Gauss} and \eqref{Codazzi}. Thus, by 
Theorem \ref{fund} there is an infinitesimal bending $\T$ of 
$f$ having $(\beta,0)$ as associated pair.\qed

\section{Variations of a product of manifolds}

In this section, we consider infinitesimal variations of 
submanifolds that are intrinsically a Riemannian product
of manifolds.
\vspace{2ex}

Let $M^n=M_1^{n_1}\times\cdots\times M_r^{n_r}$ be a 
Riemannian product of manifolds of dimensions $n_i\geq 2$, 
$1\leq i\leq r$. Given isometric immersions 
$f_i\colon M_i^{n_i}\to\R^{m_i}$, $1\leq i\leq r$, then the 
\emph{extrinsic product} $f\colon M^n\to\R^m$ of the submanifolds 
is the isometric immersion given by
$$
f(x)=(f_1(x_1),\dots,f_r(x_r)),
$$
where $x=(x_1,\dots,x_r)$ and $\R^m=\oplus_{i=1}^r\R_i^{m_i}$.
The normal space of $f$ at $x=(x_1,\dots,x_r)\in M^m$ is given by
$$
N_fM(x)=\oplus_{i=1}^rN_{f_i}M_i(x_i).
$$
where $N_{f_i}M_i(x_i)$ is the normal space of $f_i$ at 
$x_i\in M_i^{n_i}$, $1\leq i\leq r$.

Let $\iota_i^{\bar{x}}\colon M_i^{n_i}\to M^n$ denote the inclusion 
map for $\bar{x}=(\bar{x}_1\dots,\bar{x}_r)$, that is,
$$
\iota_i^{\bar{x}}(x_i)=(\bar{x}_1,\dots,x_i,\dots,\bar{x}_r),
$$
and let $\tilde{\iota}^y_i$ be the inclusion of $\R^{m_i}$ 
into $\R^m$ defined in a similar manner.
Then the second fundamental form $\a$ of $f$ at $x$ satisfies
\be\label{adapted}
\a(\iota_{i*}^x X, \iota_{j*}^x Y)=\begin{cases}
\tilde{\iota}^{f(x)}_{i*}\a_i(X,Y)\;\; &\mbox{if}\;i=j\\
0& \mbox{if}\; i\neq j
\end{cases}
\ee
for any $X\in\mathfrak{X}(M_i)$ and $Y\in\mathfrak{X}(M_j)$,
where $\a_i$ is the second fundamental form of $f_i$. 

Let $\mathcal{T}_i$ be an infinitesimal bending of 
$f_i$ in $\R^{m_i}$ for each $1\leq i\leq r$. Then,
$\mathcal{T}(x)
=\sum_{i=1}^r\tilde{\iota}_{i*}^{f(x)}\mathcal{T}_i(x_i)$ 
is an infinitesimal bending of $f$ in $\R^m$.
Let $L$ be associated to $\mathcal{T}$ and let $L_i$ be 
associated to $\mathcal{T}_i$. Then  
\be\label{Ladapted}
L\iota_{i*}^xX=\tilde{\iota}_{i*}^{f(x)}L_iX
\ee
for any $X\in\mathfrak{X}(M_i)$. If $B_i$ is associated to 
$\mathcal{T}_i$, it follows that
$$
B(\iota_{i*}^xX,\iota_{j*}^xY)
=(\nab_{\iota_{i*}^xX}L)\iota_{j*}^xY=
\begin{cases}
\tilde{\iota}_{i*}^{f(x)} B_i(X,Y)\;\;&\mbox{if}\;i=j\\
0&\mbox{if} \;i\neq j
\end{cases}
$$
for any $X\in\mathfrak{X}(M_i)$ and $Y\in \mathfrak{X}(M_j)$. 
In particular,
\be\label{betadapted}
\beta(\iota_{i*}^xX,\iota_{j*}^xY)=\begin{cases}
\tilde{\iota}_{i*}^{f(x)}\beta_i(X,Y)\;\; 
&\mbox{if}\; i=j\\
0&\mbox{if} \;i\neq j
\end{cases}
\ee
and
\be\label{eadapted}
\mathcal{E}(\iota_{i*}^xX,\tilde{\iota}_{j*}^{f(x)}\eta)=
\begin{cases}
\tilde{\iota}_{i*}^{f(x)}\mathcal{E}_i(X,\eta)\;\;&\mbox{if} \;i=j\\
0&\mbox{if} \;i\neq j
\end{cases}
\ee
for any $X\in\mathfrak{X}(M_i)$, $Y\in\mathfrak{X}(M_j)$, 
$\eta\in\Gamma(N_{f_j}M_j)$ and  $(\beta_i,\mathcal{E}_i)$ 
is the pair associated to $\T_i$, where \eqref{eadapted} 
follows  from \eqref{adapted}, \eqref{Ladapted} and the 
definition of $\Y$.\vspace{1ex}

If  $(\beta,\mathcal{E})$ is the associated pair to an  
infinitesimal bending $\mathcal{T}$ of an extrinsic product 
$f=(f_1,\dots,f_r)$  we say that $\beta$
is \emph{adapted} to the product structure if 
$$
\beta(\iota_{i*}^xX,\iota_{j*}^xY)=0
$$
for any $X\in\mathfrak{X}(M_i)$ and $Y\in\mathfrak{X}(M_j)$ 
with $i\neq j$.

\begin{proposition}\label{decomp}
Let $f\colon M^n\to\R^m$ be an extrinsic product of isometric 
immersions $f_i\colon M_i^{n_i}\to \R^{m_i}$, $n_i\geq 2$, 
$1\leq i\leq r$, that has full first normal spaces. 
If the tensor $\beta$ in the pair associated to an infinitesimal bending 
$\mathcal{T}$ of $f$ is adapted, then there exist locally 
infinitesimal bendings $\T_i$ of $f_i$, $1\leq i\leq r$, such that
$\mathcal{T}(x)=\sum_{i=1}^r\tilde{\iota}_{i*}^{f(x)}\mathcal{T}_i$.
\end{proposition}

\proof From \eqref{Gauss} we obtain
\be\label{ii}
\<\beta(\iota_{i*}^xX,\iota_{i*}^x Y),\a(\iota_{j*}^xZ,\iota_{j*}^xW)\>
+\<\a(\iota_{i*}^xX,\iota_{i*}^x Y),\beta(\iota_{j*}^xZ,\iota_{j*}^xW)\>=0
\ee
for any $X,Y\in\mathfrak{X}(M_i)$ and
$Z,W\in\mathfrak{X}(M_j)$ with $i\neq j$.

Let $\a_i(X_k,Y_k)$, $1\leq k\leq\dim N_{f_i}M_i$ with
$X_k,Y_k\in\mathfrak{X}(M_i)$ be a basis of $N_{f_i}M_i$, and set
$$
C_{ij}\a_i(X_k,Y_k)
=\beta_{N_{f_j}M}(\iota_{i*}^xX_k,\iota_{i*}^xY_k),\;j\neq i,
$$
where $\beta_{N_{f_j}M_j}$ is the component of $\beta$ in $N_{f_j}M_j$.
We claim that the linear extension to a map 
$C_{ij}\colon N_{f_i}M_i\to N_{f_j}M_j$, $i\neq j$, 
satisfies 
$$
C_{ij}\a_i(X,Y)
=\beta_{N_{f_j}M_j}(\iota_{i*}^xX,\iota_{i*}^xY)
$$
for any $X,Y\in\mathfrak{X}(M_i)$. In fact, if 
$$
\a_i(X,Y)=\sum_kc_k\a_i(X_k,Y_k)
$$ 
for $X,Y\in \mathfrak{X}(M_i)$ and 
$c_k\in\R$, $1\leq k\leq \dim N_{f_i}M_i$, we obtain 
from \eqref{ii} that
$$
\<\beta(\iota_{i*}^xX,\iota_{i*}^xY)
-\sum_kc_k\beta(\iota_{i*}^xX_k,\iota_{i*}^xY_k),
\a(\iota_{j*}^xZ,\iota_{j*}^xW)\>=0
$$
for any $Z,W\in\mathfrak{X}(M_j)$, $i\neq j$, and the
claim follows.

We obtain from \eqref{ii} that the map
$C\in\Gamma(\End(N_fM))$ defined by
$$
C\tilde{\iota}_{i*}^{f(x)}\eta_i
=\sum_{j\neq i}\tilde{\iota}_{j*}^{f(x)}C_{ij}\eta_i
$$
where $\eta_i\in\Gamma(N_{f_i}M_i)$, $1\leq i\leq r$,
is skew-symmetric. Then, we have that 
$\beta(\iota_{i*}^xX,\iota_{i*}^xY)$ decomposes orthogonally as
\be\label{betaprod}
\beta(\iota_{i*}^xX,\iota_{i*}^xY)
=\beta_{N_{f_i}M_i}(\iota_{i*}^xX,\iota_{i*}^xY)
+C\a(\iota_{i*}^xX,\iota_{i*}^xY)
\ee
for any $X,Y\in \mathfrak{X}(M_i)$.

Let $L_i\colon TM_i\to f_i^*T\R^{m_i}$ be given by
$$
L_iX=(L\iota_{i*}^xX)_{\R^{m_i}}.
$$
Since $f$ is an extrinsic product of immersions, we have 
\begin{align*}
\nab_{\iota_{j*}^xY} \tilde{\iota}_{i*}^{f(x)}L_iX
&=\nab_{\iota_{j*}^xY}(L\iota_{i*}^xX)_{\R^{m_i}}
=(\nab_{\iota_{j*}^xY}L\iota_{i*}^xX)_{\R^{m_i}}
=(B(\iota_{j*}^xY,\iota_{i*}^xX))_{\R^{m_i}}\\
&=(\Y\a(\iota_{j*}^xY,\iota_{i*}^xX)
+\beta(\iota_{j*}^xY,\iota_{i*}^xX))_{\R^{m_i}}=0
\end{align*}
for any $X\in\mathfrak{X}(M_i)$ and $Y\in\mathfrak{X}(M_j)$
with $i\neq j$, where the last steps follow using \eqref{tam} and the
assumption on $\beta$.
Thus the tensors $L_i$ are well defined on $M_i$, $1\leq i\leq r$. 
Moreover, since $B$ is symmetric, these tensors verify
$$
(\nab^i_XL_i)Y=(\nab^i_YL_i)X
$$
for any $X,Y\in\mathfrak{X}(M_i)$, where $\nab^i$ is the connection 
in $\R^{m_i}$. Thus, locally there exist vector fields 
$\T_i\in\Gamma(f_i^*T\R^{m_i})$ with $\nab^i_X\T_i=L_iX$
for any $X\in\mathfrak{X}(M_i)$, $1\leq i\leq r$.
In particular, since $L_i$ verifies \eqref{inf}, $\T_i$ is an 
infinitesimal bending of $f_i$ and,
if $\beta_i$ is associated to $\T_i$, we have
\be\label{betas}
\tilde{\iota}_{i*}^{f(x)}\beta_i(X,Y)
=\beta_{N_{f_i}M_i}(\iota_{i*}^xX,\iota_{i*}^xY)
\ee
for any $X,Y\in\mathfrak{X}(M_i)$, $1\leq i\leq r$.

Define an infinitesimal bending $\tilde{\T}$ of $f$ by
$\tilde{\T}=\sum_{i=1}^r\tilde{\iota}_{i*}^{f(x)}\mathcal{T}_i$.
We have from \eqref{betadapted}, \eqref{betaprod} and 
\eqref{betas} that $\T-\tilde{\T}$ has the associated 
tensor $\beta-\tilde{\beta}= C\a$. Since $C$ is skew-symmetric  
then the tensor $\nap C$ satisfies \eqref{anti}. Moreover, 
we have 
\begin{align*}
(\nap_X\beta-\nap_X\tilde{\beta})(Y,Z)&=(\nap_X C\a)(Y,Z)\\
&=(\nap_X C)\a(Y,Z)+C(\nap_X\a)(Y,Z)
\end{align*}
for any $X,Y,Z\in\mathfrak{X}(M)$. Using the Codazzi equation, 
we obtain
$$
(\nap_X\beta-\nap_X\tilde{\beta})(Y,Z)
-(\nap_Y\beta-\nap_Y\tilde{\beta})(X,Z)
=(\nap_X C)\a(Y,Z)-(\nap_Y C)\a(X,Z)
$$
for any $X,Y,Z\in\mathfrak{X}(M)$. Now Proposition \ref{uniqueE} 
gives $\mathcal{E}-\tilde{\mathcal{E}}=-(\nap C)$, and the proof 
follows from Proposition~\ref{trivial}.\qed

\begin{remark} {\em If for  $f\colon M^n\to\R^m$
it fails that the first normal spaces $N_1^f$ are full then
any smooth normal vector field in $N_1^\perp$ is an infinitesimal 
bending.}
\end{remark}

The \emph{$s$-nullity} $\nu_s(x)$, $1\leq s\leq m-n$, 
at $x\in M^n$ of an isometric immersion $f\colon M^n\to\R^m$ 
is defined as
$$
\nu_s(x)=\max\{\dim\mathcal{N}(\a_{U^s})(x)\colon U^s\subset N_fM(x)\}
$$
where $\a_{U^s}=\pi_{U^s}\circ\a$, $\pi_{U^s}\colon N_fM\to U^s$
is the orthogonal projection and
$$
\mathcal{N}(\a_{U^s})(x)=\{Y\in T_xM:\a_{U^s}(Y,X)=0
\;\;\mbox{for all}\;\;X\in T_xM\}.
$$

\begin{theorem}\label{snulprod}
Let $f\colon M^n\to\R^{n+p}$, $p<n$, be an extrinsic product of 
isometric immersions $f_i\colon M_i^{n_i}\to \R^{n_i+p_i}$, $n_i\geq 2$
and $1\leq i\leq r$. Assume that the \mbox{$s$-nullities} of $f$ satisfy
$\nu_s<n-s$, $1\leq s\leq p$, at any point of $M^n$. 
Then any infinitesimal bending $\T$ of $f$ is locally 
of the form 
$\T(x)=\sum_{i=1}^r\tilde{\iota}_{i*}^{f(x)}\mathcal{T}_i$, 
where $\T_i$, $1\leq i\leq r$, is an infinitesimal bending 
of $f_i$.
\end{theorem}

\proof We argue that the associated tensor $\beta$ to any 
infinitesimal bending $\T$ of $f$ has to be adapted.
Since $f$ is an extrinsic product of immersions, then we have 
$\sum_{j\neq i}\iota_{j*}^xTM_j\subset\mathcal{N}(\a_{N_{f_i}M_i})$.
Thus the assumption on the $s$-nullities yields
$\sum_{j\neq i}n_j< n-p_i$, that is,
\be\label{codim}
p_i<n_i,\;\;1\leq i\leq r.
\ee

We assume $X\in\mathfrak{X}(M_i)$, $W\in\mathfrak{X}(M_j)$ and 
$Y,Z\in\mathfrak{X}(M_k)$  with $i\neq j$. If also $k\neq i,j$, 
we obtain from \eqref{Gauss} that
\be\label{ijk}
\<\beta(\iota_{i*}^xX,\iota_{j*}^xW),
\a(\iota_{k*}^xY,\iota_{k*}^xZ)\>=0.
\ee
On the other hand, for $k=i$ it follows from \eqref{Gauss} that
\be\label{ij}
\<\beta(\iota_{i*}^xX,\iota_{j*}^xW),
\a(\iota_{i*}^xY,\iota_{i*}^xZ)\>
-\<\beta(\iota_{i*}^xY,\iota_{j*}^xW),
\a(\iota_{i*}^xX,\iota_{i*}^xZ)\>=0.
\ee
Let $\beta^i_W\colon TM_i\to N_{f_i}M_i$ be given by
$$
\beta^i_WX=\beta(\iota_{i*}^xX,\iota_{j*}^xW)_{N_{f_i}M_i}.
$$
Suppose that $\dim\Ima\beta_W^i=s>0$. Then \eqref{codim} gives
$\dim\ker\beta^i_W=n_i-s>0$. It follows from \eqref{ij} that
$$
\<\beta(\iota_{i*}^xX,\iota_{j*}^xW),\a(\iota_{i*}^xT,\iota_{i*}^xZ)\>=0
$$
for any $T$ in $\ker\beta^i_W$. This implies that $\nu_s\geq n-s$,
which contradicts our assumption and proves that 
$\beta(\iota_{i*}^xX,\iota_{j*}^xW)_{N_{f_i}M_i}=0$ for any
$X\in\mathfrak{X}(M_i)$ and $W\in\mathfrak{X}(M_j)$.
This together with \eqref{ijk} imply that
$$
\beta(\iota_{i*}^xX,\iota_{j*}^xW)=0\;\mbox{if}\;i\neq j.
$$
Thus $\beta$ is adapted, and the proof 
follows from Proposition \ref{decomp}.\vspace{2ex}\qed

The following is the main result of this section.

\begin{theorem}
Let $f\colon M^n\to \R^{n+p}$, $2p<n$, be an isometric 
immersion of a Riemannian product 
$M^n=M_1^{n_1}\times\cdots\times M_r^{n_r}$ with $n_j\geq 2$, 
$1\leq j\leq r$. Assume that the $s$-nullities of $f$ satisfy 
$\nu_s<n-2s$, $1\leq s\leq p$, at any point of $M^n$.   
Then $f$ is an extrinsic product  of isometric immersions 
$f=(f_1,\dots,f_r)$ and any infinitesimal bending $\T$ of $f$ 
is locally of the form 
$\T(x)=\sum_{i=1}^r\tilde{\iota}_{i*}^{f(x)}\mathcal{T}_i$, 
where $\T_i$ is an infinitesimal bending of 
$f_i\colon M^{n_i}_i\to \R^{m_i}$, $1\leq i\leq r$.
\end{theorem}

\proof From Theorem $5$ in \cite{DV0} or Theorem $8.14$ 
in \cite{DT} we obtain that $f$ is an extrinsic product of 
isometric immersions, and the proof follows from 
Theorem \ref{snulprod}.\qed

\section*{Acknowledgments}

This work is the result of the visit 21171/IV/19 funded by the
Fundación Séneca-Agencia de Ciencia y Tecnología de la
Región de Murcia in connection with the ``Jiménez De La Espada"
Regional Programme For Mobility, Collaboration And Knowledge Exchange.

Marcos Dajczer was partially supported by the Fundación Séneca Grant\\
21171/IV/19 (Programa Jiménez de la Espada), MICINN/FEDER project
PGC2018-097046-B-I00, and Fundación Séneca project 19901/GERM/15, 
Spain.

Miguel I. Jimenez thanks the mathematics department of the Universidad 
de Murcia for the hospitality during his visit where part of this 
work was developed.

\noindent Marcos Dajczer\\
IMPA -- Estrada Dona Castorina, 110\\
22460--320, Rio de Janeiro -- Brazil\\
e-mail: marcos@impa.br

\bigskip

\noindent Miguel Ibieta Jimenez\\
Instituto de Ciências Matemáticas e de Computação\\
Universidade de São Paulo\\
São Carlos\\
SP 13560-970-- Brazil\\
e-mail: mibieta@impa.br
\end{document}